\input amstex  
\input amsppt.sty     
\hsize 13cm
\vsize 20cm
\def\nmb#1#2{#2}         
\def\cit#1#2{\ifx#1!\cite{#2}\else#2\fi} 
\def\totoc{}             
\def\idx{}               
\def\ign#1{}             
  
\redefine\o{\circ}  
  
\define\al{\alpha}  
\define\be{\beta}  
\define\ga{\gamma}  
\define\de{\delta}

\define\ka{\kappa}  
\define\la{\lambda}  
\define\rh{\rho}  
\define\si{\sigma}  
\define\ta{\tau}

\define\ps{\psi}  
\define\om{\omega}  
  
\define\De{\Delta}  
  
\define\La{\Lambda}

\define\row#1#2#3{#1_{#2},\ldots,#1_{#3}}  
\define\rowup#1#2#3{#1^{#2},\ldots,#1^{#3}}  
\define\x{\times}  
\define\Lam#1#2#3{\La^{#1}(#2;#3)}  
\redefine\S{{\Cal S}}  
\define\A{{\Cal A}}  
\define\s#1#2{\operatorname{sign}(#1,\bold #2)}  
  
\define\Der{\operatorname{Der}}  
\define\End{\operatorname{End}}  
  
\def\today{\ifcase\month\or  
 January\or February\or March\or April\or May\or June\or  
 July\or August\or September\or October\or November\or December\fi  
 \space\number\day, \number\year}  
\topmatter  
\title $n$-ary Lie    
and Associative Algebras  \endtitle  
\author   
Peter W. Michor \\  
Alexandre M. Vinogradov \endauthor  
\affil  
Erwin Schr\"odinger International Institute of Mathematical Physics,   
Wien, Austria\\  
Universit\'a di Salerno, Italy  
\endaffil  
\address{P. Michor: Institut f\"ur Mathematik, Universit\"at Wien,  
Strudlhofgasse 4, A-1090 Wien, Austria}\endaddress  
\email{peter.michor\@esi.ac.at}\endemail  
\address{A. M. Vinogradov:   
Dip.\ Ing.\ Inf\. e Mat.\, Universit\'a di Salerno,   
Via S.\ Allende,  
84081 Baronissi, Salerno, Italy  
}\endaddress  
\email{vinograd\@ponza.dia.unisa.it}\endemail  
\date{\today}\enddate  
\dedicatory  
To Wlodek Tulczyjew, on the occasion of his 65th birthday. 
\enddedicatory 
\thanks{Supported by Project P 10037 PHY  
of `Fonds zur F\"orderung der wissenschaftlichen Forschung'} 
\endthanks  
\endkeywords  
\abstract{With the help of the multigraded Nijenhuis-- Richardson   
bracket and the multigraded Gerstenhaber bracket from \cit!{7} for   
every $n\ge 2$ we define $n$-ary associative algebras and their   
modules and also $n$-ary Lie algebras and their modules, and we give   
the relevant formulas for Hochschild and Chevalley cohomogy.}  
\endabstract  
\endtopmatter  
  
\leftheadtext{\smc Michor, Vinogradov}  
\rightheadtext{\smc $n$-ary Lie and associative algebras}  
\document  
  
\heading Table of contents \endheading  
\noindent 1. Introduction \leaders \hbox to 
1em{\hss .\hss }\hfill {\eightrm 1}\par  
\noindent 2. Review of binary algebras and bimodules 
\leaders \hbox to 1em{\hss .\hss }\hfill {\eightrm 2}\par  
\noindent 3. $n$-ary $G$-graded associative algebras and $n$-ary 
modules \leaders \hbox to 1em{\hss .\hss }\hfill {\eightrm 5}\par  
\noindent 4. Review of $G$-graded Lie algebras and 
modules\leaders \hbox to 1em{\hss .\hss }\hfill {\eightrm 7}\par  
\noindent 5. $n$-ary $G$-graded Lie algebras and their modules 
\leaders \hbox to 1em{\hss .\hss }\hfill {\eightrm 10}\par  
\noindent 6. Relations between $n$-ary algebras and Lie 
algebras\leaders \hbox to 1em{\hss .\hss }\hfill {\eightrm 12}\par  
\noindent 7. Hochschild operations and non commutative differential 
calculus\leaders \hbox to 1em{\hss .\hss }\hfill {\eightrm 12}\par  
\noindent 8. Remarks on Filipov's $n$-ary Lie algebras 
\leaders \hbox to 1em{\hss .\hss }\hfill {\eightrm 15}\par  
\noindent 9. Dynamical aspects\leaders \hbox to 
1em{\hss .\hss }\hfill {\eightrm 17}\par
 
\head\totoc\nmb0{1}. Introduction \endhead  
  
In 1985 V.\ Filipov \cit!{3} proposed a generalization of the   
concept of a Lie algebra by replacing the binary operation by an   
$n$-ary one. He defined an $n$-ary Lie algebra structure on a vector   
space $V$ as an operation which associates with each $n$-tuple   
$(u_1,\dots,u_n)$ of elements in $V$ another element   
$[u_1,\dots,u_n]$ which is $n$-linear, skew symmetric, and satisfies   
the $n$-Jacobi identity:  
$$  
[u_1,\dots,u_{n-1},[v_1,\dots,v_n]] =   
\sum [v_1,\dots,v_{i-1}[u_1,\dots,u_{n-1},v_i],\dots,v_n].  
\tag1$$   
Apparently Filippov was motivated by the fact that with this   
definition one can delelop a meaningful structure theory, in   
accordance with the aim of Malcev's school: To look for algebraic   
structures that manifest good properties.  
  
On the other hand, in 1973 Y.\ Nambu \cit!{13}  
proposed an $n$-ary generalization of Hamiltonian dynamics by means   
of the $n$-ary `Poisson bracket'  
$$  
\{f_1,\dots,f_n\} = \det\left(\frac{\partial f_i}{\partial x_j}\right).  
\tag2$$  
Apparently he looked for a simple model which explains the   
unseparability of quarks. Much later, in the early 90's, it was   
noticed by M.\ Flato, C.\ Fronsdal, and others, that the $n$-bracket   
\thetag 2 satisfies \thetag 1. On this basis L.\ Takhtajan \cit!{17}  
developed sytematically the foundations of of the theory of   
$n$-Poisson or Nambu-Poisson manifolds. It seems that the work of   
Filippov was unknown then; in particular Takhtajan reproduces some   
results from \cit!{3} without refereing to it.   
  
Recently Alekseevsky and Guha \cit!{1} and later Marmo, Vilasi, and   
Vinogradov \cit!{9} proved that $n$-Poisson structures of the kind   
above are extremely rigid: Locally they are given by $n$ commuting vector 
fields   of rank $n$, if $n>2$; in other words, $n$-Poisson structures are   
locally given by \thetag2. This rigidity suggests that one should   
look for alternative $n$-ary analogs of the concept of a Lie algebra.   
One of them is proposed below in this paper. It is based on the   
completely skew symmetrized version of Filippov's Jacobi identity   
\thetag2. It is shown in \cit!{20} that this approach leads to richer   
and more diverse structures which seem to be more useful for purposes   
of dynamics. In fact, we were lead in 1990-92 to the constructions of   
this paper by some expectations about $n$-body mechanics and the   
naturality of the machinary developed in \cit!{7}. So, our motives  
were quite different from that by Filippov, Nambu and Takhtajian. 
This paper is essentailly based on our   unpublished notes from 
 1990-92. In view 
of the recent developments we   decided to publish them now. In this paper we 
consider $G$-graded   $n$-ary generalizations of the concept of associative 
algebras, of   Lie algebras, their modules, and their cohomologies; all this is   
produced by the algebraic machinery of \cit!{7}. Related (but not   
graded) concepts are discussed in \cit!{4} in terms of operads and   
their Koszul duality.  
The recent preprints \cit!{2} and \cit!{5} propose dynamical models   
which correspond to the not graded case with even $n$  in our construction.  
   
\heading\totoc\nmb0{2}. Review of binary algebras and bimodules \endheading  
  
In this section we review the results from the paper \cit!{7} in  
a slightly different point of view.  
  
\subheading{\nmb.{2.1}. Conventions and definitions} By a {\it  
grading group} we mean a commutative group $(G,+)$ together with a   
$\Bbb Z$-bilinear symmetric mapping (bicharacter)  
$\langle \quad,\quad\rangle:G\x G\to \Bbb Z_2:=\Bbb Z/2\Bbb Z$.   
Elements of $G$ will be called degrees, or $G$-degrees if more   
precision is necessary. A standard example of a grading group is   
$\Bbb Z^m$ with   
$\langle x,y\rangle  = \sum_{i=1}^m x^iy^i(\mod 2)$.  
If $G$ is a grading group we will consider the grading group   
$\Bbb Z\x G$ with   
$\langle (k,x),(l,y)\rangle = kl (\mod2)+\langle x,y\rangle$.  
  
A {\it $G$-graded vector space} is just a direct sum $V =  
\bigoplus_{x \in G}V^x$, where the elements of $V^x$ are  
said to be homogeneous of $G$-degree $x$. We assume  
that vector spaces are defined  
over a field $\Bbb K$ of characteristic 0.  
In the following $X$,  
$Y$, etc will always denote homogeneous elements of some  
$G$-graded vector space of $G$-degrees $x$, $y$, etc.  
  
By an $G$-graded algebra $\A  = \bigoplus_{x \in  
G}\A ^x$ we mean an $G$-graded vector space which is also  
a $\Bbb K$ algebra such that $\A^x\cdot \A^y \subseteq \A^{x+y}$.  
\roster \item The $G$-graded algebra $(\A,\cdot)$ is said to be {\it  
$G$-graded commutative} if for homogeneous elements $X$, $Y \in  
\A $ of $G$-degree $x$, $y$, respectively,we have $X\cdot Y =  
(-1)^{\langle x,y\rangle }Y\cdot X$.  
\item If $X\cdot Y = -(-1)^{\langle x,y\rangle }Y\cdot X$ holds it is   
     called {\it $G$-graded anticommutative}.  
\item By an {\it $G$-graded Lie algebra} we mean a $G$-graded  
     anticommutative algebra $(\Cal E, [\quad,\quad])$ for which the {\it  
     $G$-graded Jacobi identity} holds:   
     $$[X,[Y,Z]] = [[X,Y],Z] + (-1)^{\langle x,y\rangle }[Y,[X,Z]]$$  
\endroster  
Obviously the space   
$\End(V) = \bigoplus_{\de \in G}\End^{\de}(V)$  
of all endomorphisms of a $G$-graded vector space  
$V$ is a $G$-graded algebra under composition, where  
$\End^\de(V)$ is the space of linear endomorphisms $D$ of $V$ of  
$G$-degree $\de$, i.e. $D(V^x) \subseteq V^{x+\de}$. Clearly  
$\End(V)$ is a $G$-graded Lie algebra under the $G$-graded  
commutator  
$$[D_1,D_2] := D_1\circ D_2 - (-1)^{\langle \de_1,\de_2\rangle } D_2\circ D_1.\tag4$$  
  
If $\A $ is a $G$-graded algebra, an endomorphism   
$D:\A \to \A$ of $G$-degree $\de$ is called a {\it $G$-graded  
derivation}, if for $X$, $Y \in \A $ we have   
$$D(X\cdot Y) = D(X) \cdot Y + (-1)^{\langle \delta,x\rangle }X\cdot D(Y).\tag5$$   
Let us write $\Der^\de (\A )$ for the space of all $G$-graded  
derivations of degree $\de$ of the algebra $\A $, and we put   
$$\Der(\A ) = \bigoplus_{\de \in G}\Der^\de (\A ).\tag5$$  
  
The following lemma is standard:  
\proclaim{Lemma} If $\A $ is an $G$-graded algebra, then the space  
$\Der(\A )$ of $G$-graded derivations is an $G$-graded Lie   
algebra under the $G$-graded commutator.  
\endproclaim  
  
\subheading{\nmb.{2.2} Graded associative algebras}   
Let $V=\bigoplus_{x\in G}V^x$ be   
an $G$-graded vector space. We define   
$$M(V):=\bigoplus_{(k,\ka)\in\Bbb Z\x G}M^{(k,\ka)}(V),$$  
where $M^{(k,\ka)}(V)$ is the space of all $k+1$-linear mappings  
$K:V\x\dots\x V\to V$ such that   
$K(V^{x_0}\x\dots\x V^{x_k})\subseteq V^{x_0+\dots+x_k+\ka}$.  
We call $k$ the \idx{\it form degree} and $\ka$ the \idx{\it weight   
degree} of $K$.  
We define for $K_i\in M^{(k_i,\ka_i)}(V)$ and $X_j\in V^{x_j}$  
$$\gather(j(K_1)K_2)(\row X0{k_1+k_2}) := \\  
=\sum_{i=0}^{k_2} (-1)^{k_1i+\langle \ka_1,\ka_2+x_0+\dots+x_{i-1}\rangle}  
	K_2(X_0,\dots,K_1(\row Xi{i+k_1}),\dots,X_{k_1+k_2}),\\  
[K_1,K_2]^\De=j(K_1)K_2-(-1)^{k_1k_2+\langle \ka_1,\ka_2\rangle}j(K_2)K_1.  
\endgather$$  
  
\proclaim{Theorem} Let $V$ be an $G$-graded vector  
space. Then we have:  
\roster  
\item $(M(V),[\quad,\quad]^\De)$ is a $(\Bbb Z\x G)$-graded Lie algebra.  
\item If $\mu\in M^{(1,0)}(V)$, so $\mu:V\x V\to V$ is bilinear   
     of weight $0\in G$, then $\mu$ is an associative $G$-graded   
     multiplication if and only if $j(\mu)\mu=0$.  
\item If $\nu\in M^{(1,n)}(V)$, so $\nu:V\x V\to V$ is bilinear   
     of weight $n\in G$, then $j(\nu)\nu=0$ is equivalent to  
$$ 
\nu(\nu(X_0,X_1),X_2) -(-1)^{\langle n,n\rangle}\nu(X_0,\nu(X_1,X_2)) =0  
$$ 
     which is the natural notion of an associative multiplication of  
     weigth $n\in G$. 
\endroster  
\endproclaim  
  
\demo{Proof} The first assertion is from \cit!{7}.  
The second and third assertion follows by writing out the definitions.  
\qed\enddemo  
 
In \cit!{7} the formulation was as follows:  
$\mu\in M^{(1,0)}(V)$ is an associative $G$-graded algebra structure if and   
only if $[\mu,\mu]^{\De}=2j(\mu)\mu=0$.  
For $\nu\in M^{(1,n)}(V)$ we have  
$[\nu,\nu]^{\De}=(1+(-1)^{\langle n,n\rangle})j(\nu)\nu$.

\subheading{\nmb.{2.3}. Multigraded bimodules}  
Let $V$ and $W$ be $G$-graded vector spaces and $\mu:V\x V \to V$ a   
$G$-graded algebra structure. A {\it $G$-graded bimodule}   
$\Cal M=(W,\la,\rh)$ over $\A =(V,\mu)$  
is given by $\la ,\rh : V\to \End(W)$ of weight $0$ such that  
$$\align  
j(\mu)\mu &= 0 \quad \text{ so } \A \text{ is associative } \tag1 \\  
\la (\mu (X_1,X_2)) &= \la (X_1)\circ \la (X_2)  \tag2 \\  
\rho (\mu (X_1,X_2)) &= (-1)^{\langle x_1,x_2\rangle}  
     \rho (X_2)\circ \rho (X_1)  \tag3 \\  
\la (X_1)\circ \rho (X_2) &= (-1)^{\langle x_1,x_2\rangle }  
     \rho (X_2)\circ \la (X_1) \tag4   
\endalign$$  
where $X_i \in V^{x_i}$ and $\circ$ denotes the composition in $\End (W)$.  
  
\proclaim{\nmb.{2.4}. Theorem}  
Let $E$ be the $(\Bbb Z\x G)$-graded vector space defined by  
$$E^{(k,*)}=\left\{\alignedat2 &V &\qquad &\text{if } k=0 \\   
                    &W &&\text{if } k=1 \\   
                    &0 && \text{otherwise.}   
               \endalignedat  
\right. $$  
Then $P\in M^{(1,0)}(E)$ defines a bimodule structure on $W$   
if and only if $j(P)P=0$.  
\endproclaim  
  
\demo{Proof}  
We define  
$$\align  
\mu (X_1,X_2)&:=P(X_1,X_2) \\  
\la (X)Y&:=P(X,Y) \\  
\rho(X)Y&:=(-1)^{\langle x,y\rangle }P(Y,X)  
\endalign$$  
where we suppose the $X_i$'s $\in V$ and $Y \in W$ to be embedded in $E$.  
Then if $Z_i \in E$ is arbitrary we get  
$$(j(P)P)(Z_0,Z_1,Z_2) = P((Z_0,Z_1),Z_2) - P(Z_0,(Z_1,Z_2)). $$  
Now specify $Z_i\in V$ resp. $W$ to get eight independent equations. Four  
of them vanish identically because of their degree of homogeneity, the   
others recover the defining equations for the $G$-graded bimodules.  
\qed\enddemo  
  
\proclaim{\nmb.{2.5} Corollary}  
In the above situation we have the following decomposition of the   
$(\Bbb Z^2\x G)$-graded space $M(E)$ :  
$$M^{(k,q,*)}(E) = \left\{  
     \alignedat2 &0 &\quad &\text{for } q>1 \\  
          &L^{(k+1,*)}(V,W) &&\text{for } q=1 \\  
          &M^{(k,*)}(V)\bigoplus^{k+1}(L^{(k,*)}(V,\End(W))  
          &&\text{for } q=0   
     \endalignedat  
\right.$$  
where $L^{(k,*)}(V,W)$ denotes the space of $k$-linear mappings   
$V\x\dots\x V\to W$.   
If $P$ is as above, then  
$P=\mu +\la +\rho$ corresponds exactly to this decomposition.  
\qed\endproclaim  
  
\subheading{\nmb.{2.6}. Hochschild cohomology and multiplicative   
structures}  
Let $V$,$W$ and $P$ be as in {Theorem \nmb!{2.4}} and let $\nu : W\x W \to W$ be  
a $G$-graded algebra structure, so $\nu \in 
M^{(1,-1,0)}(E)$.  Then for $C_i\in L^{(k_i,c_i)}(V,W)$ we define  
$$C_1\bullet C_2:=[C_1,[C_2,\nu]^{\De}]^{\De} = \pm \nu(C_1,C_2).$$  
Since $[C_1,C_2]^{\De}=0$ it follows that $(L(V,W),\bullet)$ is   
$(\Bbb Z\x G)$-graded commutative.  
  
\proclaim{Theorem}  
 
1. The mapping $[P,\quad]^{\De} : M(E)\to M(E)$ is a differential. Its 
restriction $\de_P$ to $L(V,W)$  is a generalization of the 
Hochschild coboundary operator to the $G$-graded case: If  $C \in 
L^{(k,c)}(V,W)$, then we have for $X_i \in V^{x_i}$  $$\multline   
(\de_P C)(\row X0k) =  \la (X_0)C(\row X1k) \\  
-\sum_{i=0}^{k-1}(-1)^i C(X_0,\dots,\mu (X_i,X_{i+1}),\dots,X_k) \\  
+(-1)^{k+1+\langle x_0+\cdots+x_{k-1}+c,x_k\rangle}\rho (X_k)C(\row X0{k-1})  
\endmultline$$  
The corresponding $(\Bbb Z\x G)$-graded cohomology will be denoted by   
$H(\A,\Cal M)$.   
  
2. If $[P,\nu]^{\De}=0$, then $\de_P$ is a derivation of $L(V,W)$ of   
$(\Bbb Z\x G)$-degree $(1,0)$. In this case the product $\bullet$ carries over  
to a $(\Bbb Z\x G)$-graded (cup) product on $H(\A,\Cal M)$.  
\endproclaim  
  
\heading\totoc\nmb0{3}. $n$-ary $G$-graded associative algebras   
 and $n$-ary modules \endheading  
  
\subheading{\nmb.{3.1}. Definition} Let $V$ be a $G$-graded vector   
space. Let $\mu\in M^{(n-1,0)}(V)$, so $\mu: V^{\otimes n}\to V$   
is $n$-linear of weight $0\in G$.   
  
We call $\mu$ an   
\idx{\it $n$-ary associative $G$-graded multiplication} of weigth  
$0\in G$ if  $j(\mu)\mu=0\in M^{(2n-2,0)}(V)$.  
  
\subheading{Remark} We are forced to use $j(\mu)\mu=0$ instead of  
$[\mu,\mu]^{\De}=0$ since the latter condition is automatically 
satisfied   for odd $n$.  
  
\subheading{\nmb.{3.2}. Example} If $V$ is 0-graded, then a ternary   
associative multiplication $\mu: V\x V\x V\to V$ satisfies  
$$\multline  
(j(\mu)\mu)(\row X05)  
= \mu(\mu(X_0,X_1,X_2),X_3,X_4)+\\  
+ \mu(X_0,\mu(X_1,X_2,X_3),X_4)  
+ \mu(X_0,X_1,\mu(X_2,X_3,X_4)) =0.  
\endmultline$$  
  
\subheading{\nmb.{3.3}. Definition} Let $V$ and $W$ be $G$-graded   
vector spaces. We consider the $(\Bbb Z\x G)$-graded vector space   
$E$ defined by  
$$E^{(k,*)}=\left\{\alignedat2 &V &\qquad &\text{if } k=0 \\   
                    &W &&\text{if } k=1 \\   
                    &0 && \text{otherwise.}   
               \endalignedat  
\right. $$  
Then $P\in M^{(n-1,0,0)}(E)$ is called an   
\idx{\it $n$-ary $G$-graded module structure} on $W$ over an  
$n$-ary algebra structure on $V$   
if  $j(P)P=0$.  
Let us denote the resulting $n$-ary algebra by $\A$,  
and the $n$-ary module by $\Cal W$.  
  
The mapping $P$ is the sum of partial mappings  
$$\align  
&\mu=P: V\x\dots\x V \to V\qquad\text{ the $n$-ary algebra structure }\\  
&P: W\x V\x \dots \x V \to W \quad\text{ the rightmost $n$-ary   
     module structure} \\  
&P: V\x W\x V\x \dots \x V \to W \\  
&\quad\ldots \\  
&P: V\x\dots\x V\x W\x V \to W \\  
&P: V\x\dots\x V\x W \to W\quad\text{ the leftmost $n$-ary module   
     structure}  
\endalign$$  
This decomposition of $P$ corresponds exactly to the last line in the   
decomposition of $M^{(n-1,0,*)}$ of \nmb!{2.5}. 
 
The above definition is easily generalized  by changing the form degree
of $W$ or/and by augmenting the number of $W$'s. For simplicity we don't
discuss this possibility here.
  
\subheading{\nmb.{3.4}. Example} If $V$ and $W$ are 0-graded then   
a ternary module satisfies the following conditions besides the one   
from \nmb!{3.2} describing the ternary algebra structure on $V$:  
$$\align  
&P(P(w_0,v_1,v_2),v_3,v_4)  +P(w_0,\mu(v_1,v_2,v_3),v_4)+P(w_0,v_1,\mu(v_2,v_3,v_4))=0\\  
&P(P(v_0,w_1,v_2),v_3,v_4)  +P(v_0,P(w_1,v_2,v_3),v_4)  +P(v_0,w_1,\mu(v_2,v_3,v_4))=0\\  
&P(P(v_0,v_1,w_2),v_3,v_4)  +P(v_0,P(v_1,w_2,v_3),v_4)  +P(v_0,v_1,P(w_2,v_3,v_4))=0\\  
&P(\mu(v_0,v_1,v_2),w_3,v_4)+P(v_0,P(v_1,v_2,w_3),v_4)  +P(v_0,v_1,P(v_2,w_3,v_4))=0\\  
&P(\mu(v_0,v_1,v_2),v_3,w_4)+P(v_0,\mu(v_1,v_2,v_3),w_4)+P(v_0,v_1,P(v_2,v_3,w_4))=0  
\endalign$$  
  
\subheading{\nmb.{3.5}. Hochschild cohomology for even $n$}  
Let $V$ and $W$ be $G$-graded vector spaces,   
and let $P\in M^{(n-1,0,0)}(E)$ be an $n$-ary   
module structure on $W$ over an $n$-ary $G$-graded algebra structure   
on $V$ as in definition \nmb!{3.3}.  
  
\proclaim{Theorem} Let $n=2k$ be even. Then we have:  
  
The mapping $[P,\quad]^{\De} : M(E)\to M(E)$ is a differential.    
 Its restriction $\de_P$ to $L(V,W)$ is called 
the Hochschild coboundary operator. For a   
cochain $C\in M^{(k,1,c)}= L^{(k+1,c)}(V,W)$ and with $p=n-1$  
we have for $X_i \in V^{x_i}$  
$$\multline   
(\de_P C)(\row X0{k+p}) =   
\sum_{i=0}^{k}(-1)^{pi} C(X_0\dots,P(X_i,\dots,X_{i+p}),\dots,X_{k+p}) \\  
-\sum_{j=0}^p(-1)^{k(j+p)+\langle x_0+\cdots+x_{j-1},c\rangle}  
P(X_0,\dots,C(X_j,\dots,X_{j+k}),\dots,X_{k+p}).  
\endmultline$$  
The corresponding $(\Bbb Z\x G)$-graded cohomology will be denoted by   
$H(\A,\Cal M)$.   
\endproclaim  
  
\demo{Proof}  
We have by the $(\Bbb Z^2\x G)$-graded Jacobi identity  
$$[P,[P,Q]^{\De}]^{\De}= [[P,P]^{\De},Q]^{\De} +   
     (-1)^{(n-1)^2}[P,[P,Q]^{\De}]^{\De}$$  
which implies that $[P,\quad]^{\De}$ is a differential since $n-1$ is   
odd and $[P,P]^{\De}= j(P)P -(-1)^{(n-1)^2}j(P)P = 2j(P)P =0$.  
The rest follows from a computation.  
\qed\enddemo  
 
\subhead\nmb.{3.6}. Remark \endsubhead  
We get an easy extension of the Hochschild coboundary operator for  
$n$-ary algebra structures for odd $n$ if we choose the weigth  
accordingly.  
Let $P\in M^{(n-1,0,p)}(E)$ be an $n$-ary   
module structure of weight $p$ on $W$ over an $n$-ary $G$-graded  
algebra structure of weight $p$ on $V$, similarly as in definition  
\nmb!{3.3}: We require that $j(P)P=0$.  
Let us suppose that $\|(n-1,0,p)\|^2=(n-1)^2 + \langle p,p\rangle$ is  
odd. 
Then by \nmb!{2.2} we have  
$$\align 
[P,P]^\De &= \Bigl(1-(-1)^{(n-1)^2+\langle p,p \rangle}\Bigr) j(P)P =  
2 j(P)P = 0,\\ 
[P,[P,Q]^{\De}]^{\De}&= [[P,P]^{\De},Q]^{\De} +   
     (-1)^{(n-1)^2+\langle p,p\rangle}[P,[P,Q]^{\De}]^{\De} = 0, 
\endalign$$ 
so that we get a differential.  
A dual version of this can be seen in \nmb!{7.2}.(3) below.  
  
\subhead\nmb.{3.7}. Ideals \endsubhead  
Let $(V,\mu)$ be an $n$-ary $G$-graded associative algebra. An ideal   
$I$ in $(V,\mu)$ is a linear subspace $I\subset V$ such that  
$\mu(X_1,\dots,X_n)\in I$ whenever one of the $X_i\in I$. Then   
$\mu$ factors to an $n$-ary associative multiplication on the   
quotient space $V/I$. This quotient space is again $G$-graded, if  
$I$ is a $G$-graded subspace in the sense that   
$I=\bigoplus_{x\in G} (I\cap V^x)$.   
  
Of course any ideal $I$ is an $n$-ary module over $(V,\mu)$ which is   
$G$-graded if and only if $I$ is $G$-graded.   
Conversely, any $n$-ary module $W$ over $(V,\mu)$ is an ideal in the   
$n$-ary algebra $V\oplus W= E$ with the multiplication $P$ from   
\nmb!{3.3}. Here $P(X_1,\dots,X_n)=0$ if any two elements $X_i$ lie   
in $W$, so that $E$ may be regarded as an $G$-graded or as a   
$(\Bbb Z\x G)$-graded algebra. It could be called also the   
\idx{\it semidirect product} of $V$ and $W$.  
  
\subhead\nmb.{3.8}. Homomorphisms \endsubhead  
A linear mapping $f:V\to W$ of degree 0 between two $G$-graded   
algebras $(V,\mu)$ and $(W,\nu)$ is called a \idx{\it homomorphism of   
$G$-graded algebras} if it is compatible with the two $n$-ary   
multiplications:  
$$f(\mu(X_1,\dots,X_n))=\nu(f(X_1),\dots,f(X_n))$$  
Then the kernel of $f$ is an $n$-ary ideal in $(V,\mu)$ and the image   
of $f$ is an $n$-ary subalgebra of $(W,\nu)$ which is isomorphic to   
$V/\operatorname{ker}(f)$.  
  
Similarly we can define the notion of an $n$-ary $V$-module   
homomorphism  between two $V$-modules $W_0$ and $W_1$. Then the   
category of all ($G$-graded) $n$-ary $V$-modules and of their   
homomorphisms is an abelian category. We did not investigate the   
relation to the embedding theorem of Freyd and Mitchell.  
  
\heading\totoc\nmb0{4}. Review of $G$-graded Lie algebras and   
modules\endheading  
  
In this section we sketch the theory from \cit!{7} for $G$-graded   
Lie algebras from a slightly different angle. In this section section   
we need that the ground field $\Bbb K$ has characteristic 0.  
  
\subheading{\nmb.{4.1}. Multigraded signs of permutations}   
Let $\bold x = (\row x1k) \in G^k$ be a multi index of $G$-degrees  
$x_i \in G$ and let $\si \in \S_k$ be a  
permutation of $k$ symbols. Then we define the {\it $G$-graded  
sign} $\s \si x$ as follows:   
For a transposition $\si =(i,i+1)$ we  
put $\s \si x = -(-1)^{\langle x_i,x_{i+1}\rangle }$;   
it can be checked by combinatorics that this gives a well defined  
mapping $\s {\quad}x:\Cal S_k \to  \{-1,+1\}$.   
  
Let us write $\si x = (\row x{\si1}{\si k})$, then we have the following   
  
\proclaim{Lemma} $\s {\si\o\ta}x = \s \si x .\s \ta{\si \bold x}$.\qed  
\endproclaim  
  
\subheading{\nmb.{4.2} Multigraded Nijenhuis-Richardson algebra}  
We define the {\it $G$-graded alternator} $\al : M(V)\to M(V)$ by  
$$(\al K)(X_0,\dots,X_k) = \dsize\frac1{(k+1)!} \sum_{\si \in \S_{k+1}}  
\s \si x K(X_{\si 0},\dots,X_{\si k})\tag1$$  
for $K \in M^{(k,*)}(V)$ and $X_i \in V^{x_i}$. By lemma \nmb!{4.1} we have  
$\al^2 = \al$ so $\al$ is a projection on $M(V)$, homogeneous of   
$(\Bbb Z\x G)$-degree 0, and we set   
$$  
A(V)=\bigoplus_{(k,\ka)\in\Bbb Z\x G}A^{(k,\ka)}(V)  
=\bigoplus_{(k,\ka)\in \Bbb Z\x G}\al (M^{(k,\ka)}(V)).  
$$  
A long but straightforward computation shows that for   
$K_i \in M^{(k_i,\ka_i)}(V)$   
$$\al (j(\al K_1)\al K_2) = \al (j(K_1)K_2),$$  
so the following operator and bracket is well defined:  
$$\align  
i(K_1)K_2 :&= \dsize\frac{(k_1+k_2+1)!}{(k_1+1)!(k_2+1)!}  
     \al(j(K_1)K_2)\\  
[K_1,K_2]^{\wedge} &=\dsize\frac{(k_1+k_2+1)!}{(k_1+1)!(k_2+1)!}  
                                               \al ([K_1,K_2]^{\De}) \\  
&=i(K_1)K_2 -(-1)^{\langle (k_1\kappa_1),(k_2,\kappa_2)\rangle }i(K_2)K_1  
\endalign$$  
The combinatorial factor is explained in \cit!{7},~3.4.  
  
\proclaim{\nmb.{4.3}. Theorem}  
1. If $K_i$ are as above, then  
$$\multline  
(i(K_1)K_2)(\row X0{k_1+k_2}) = \\  
=\frac1{(k_1+1)!k_2!}  
     \sum_{\si \in \S_{k_1+k_2+1}}\s \si x   
     (-1)^{\langle \ka_1,\ka_2\rangle }\cdot\\  
\cdot K_2((K_1(X_{\si 0},\dots,X_{\si k_1}),\dots,X_{\si(k_1+k_2)}).  
\endmultline$$  
  
2. $(A(V),[\quad,\quad]^{\wedge})$ is a $(\Bbb Z\x G)$-graded Lie algebra.  
  
3. If $\mu\in A^{(1,0)}(V)$, so $\mu:V\x V\to V$ is bilinear   
$G$-graded anticommutative mapping of weight $0 \in G$, then   
$i(\mu)\mu=0$ if and only if $(V,\mu)$ is a $G$-graded Lie algebra.  
\endproclaim  
  
\demo{Proof}  
For 1 and 2 see \cit!{7}.   
  
3. Let $\mu \in A^{(1,0)}(V)$, then from 1 we see that  
$$(i(\mu)\mu)(X_0,X_1,X_2)   
     =\tfrac1{2!}\sum_{\si \in \Cal S_3}\s \si x\cdot \mu(\mu(X_{\si0},  
       X_{\si1}),X_{\si2}))$$  
which is equivalent to the $G$-graded Jacobi expression of $(V,\mu).$  
\qed\enddemo  
  
$(A(V),[\quad,\quad]^{\wedge})$ is called the {\it $(\Bbb Z\x G)$-graded   
Nijenhuis-Richardson algebra}, since $A(V)$ coincides for $G=0$   
with $Alt(V)$ of \cit!{14}.  
  
\proclaim{\nmb.{4.4}. Theorem} Let $V$ and $W$ be $G$-graded vector   
spaces. Let $E$ be the $(\Bbb Z\x G)$-graded vector space defined by  
$$E^{(k,*)}=\left\{\alignedat2 &V &\qquad &\text{if } k=0 \\   
                    &W &&\text{if } k=1 \\   
                    &0 && \text{otherwise.}   
               \endalignedat  
\right. $$  
Let $P\in A^{(1,0,0)}(E)$ then $i(P)P=0$ if and  
only if  
$$i(\mu)\mu=0\tag a$$    
so $(V,\mu)=\frak g$ is a $G$-graded Lie algebra, and   
$$\rh (\mu (X_1,X_2))Y = [\rh (X_1),\rh (X_2)]Y \tag b$$  
where $\mu (X_1,X_2)=P(X_1,X_2)\in V$ and $\rh (X)Y=P(X,Y)\in W$  
for $X$, $X_i\in V$ and $Y\in W$, and where $[\quad,\quad]$  
denotes the $G$-graded commutator in $\End (W)$. So $i(P)P=0$  
is by definition equivalent to the fact that $\Cal M :=(W,\rh)$ is a   
$G$-graded Lie-$\frak g$ module.  
  
If $P$ is as above the mapping   
$\partial_P := [P,\quad]^{\wedge}: A(E) \to A(E)$ is a differential   
and its restriction to   
$$\bigoplus_{k\in \Bbb Z}  
\La^{(k,*)}(\frak g,\Cal M):= \bigoplus_{k\in \Bbb Z}A^{(k,1,*)}(E)$$  
generalizes the  Chevalley-Eilenberg coboundary operator to the   
$G$-graded case:  
$$\align  
(\partial_PC)(\row X0k)   
     &= \sum_{i=0}^k(-1)^{\al_i(\bold x)+\langle x_i,c\rangle }  
     \rh (X_i)C(X_0,\dots,\widehat{X_i},\dots,X_k) \\  
&+\sum_{i<j}(-1)^{\al_{ij}(\bold x)}  
     C(\mu (X_i,X_j),\dots,\widehat{X_i},\dots,\widehat{X_j},\dotsc)  
\endalign$$  
where   
$$\left\{  
\aligned \al_i(\bold x)&=\langle x_i,x_0+\dots+x_{i-1}\rangle +i \\  
\al_{ij}(\bold x)&=\al_i(\bold x)+\al_i(\bold x)+\langle x_i,x_j\rangle   
\endaligned\right.$$  
We denote the corresponding $(\Bbb Z\x G)$-graded cohomology space by   
$H(\frak g,\Cal M)$.   
  
If $\nu : W\x W\to W$ is $G$-graded symmetric (so   
$\nu \in A^{(1,-1,*)}(E)$) and $[P,\nu]^{\wedge}=0$ then $\partial_P$ acts  
as derivation of $G$-degree $(1,0)$ on the $(\Bbb Z\x G)$-graded commutative  
algebra $(\La (\frak g,\Cal M),\bullet)$, where  
$$C_1\bullet C_2:=[C_1,[C_2,\nu]^{\wedge}]^{\wedge}  
     \quad C_i\in \La^{(k_i,c_i)}(\frak g,\Cal M).$$  
In this situation the product $\bullet$ carries over   
to a $(\Bbb Z\x G)$-graded symmetric  
(cup) product on $H(\frak g,\Cal M)$.  
\endproclaim  
  
\demo{Proof}  
Apply the $G$-graded alternator $\al$ to the results of \nmb!{2.3},   
\nmb!{2.4}, \nmb!{2.5}, and \nmb!{2.6}.  
\qed  
\enddemo  
  
\heading\totoc\nmb0{5}. $n$-ary $G$-graded Lie algebras   
and their modules \endheading  
  
\subheading{\nmb.{5.1}. Definition} Let $V$ be a $G$-graded  vector   
space. Let $\mu\in A^{(n-1,0)}(V)$, so $\mu: V^n\to V$ is a   
$G$-graded skew symmetric $n$-linear mapping.   
  
We call $\mu$ an \idx{\it $n$-ary $G$-graded Lie algebra structure} 
on $V$ if   
 $i(\mu)\mu = 0$.  
  
\subheading{\nmb.{5.2}. Example} If $V$ is 0-graded, then a ternary   
Lie algebra structure on $V$ is a skew symmetric trilinear mapping  
$\mu: V\x V\x V\to V$ satisfying  
$$\align  
0&= (i(\mu)\mu)(\row X04)   
= \frac1{3!\,2!}\sum_{\si\in\Cal S_3}   
     \operatorname{sign}(\si)\,   
     \mu(\mu(X_{\si0},X_{\si1},X_{\si2}),X_{\si3},X_{\si4}) \\   
&=     + \mu(\mu(X_0,X_1,X_2),X_3,X_4)  
       - \mu(\mu(X_0,X_1,X_3),X_2,X_4)\\  
&\quad + \mu(\mu(X_0,X_1,X_4),X_2,X_3)  
       + \mu(\mu(X_0,X_2,X_3),X_1,X_4)\\  
&\quad - \mu(\mu(X_0,X_2,X_4),X_1,X_3)  
       + \mu(\mu(X_0,X_3,X_4),X_1,X_2)\\  
&\quad - \mu(\mu(X_1,X_2,X_3),X_0,X_4)  
       + \mu(\mu(X_1,X_2,X_4),X_0,X_3)\\  
&\quad - \mu(\mu(X_1,X_3,X_4),X_0,X_2)  
       + \mu(\mu(X_2,X_3,X_4),X_0,X_1)    
\endalign$$

\subheading{\nmb.{5.3}. Definition}  
Let $V$ and $W$ be $G$-graded   
vector spaces. We consider the $(\Bbb Z\x G)$-graded vector space   
$E$ defined by  
$$E^{(k,*)}=\left\{\alignedat2 &V &\qquad &\text{if } k=0 \\   
                    &W &&\text{if } k=1 \\   
                    &0 && \text{otherwise.}   
               \endalignedat  
\right. $$  
Then $P\in A^{(n-1,0,0)}(E)$ is called an   
\idx{\it $n$-ary $G$-graded Lie module structure} on $W$ over an  
$n$-ary Lie algebra structure on $V$   
if  $i(P)P=0$.  
Let us denote the resulting $n$-ary Lie algebra by $\frak g$,  
and the $n$-ary module by $\Cal W$.  
  
Ordering by degree and using the $G$-graded skew symmetry   
we see that $P$ is now the sum of only two   
partial $n$-linear mappings   
$$\align  
&\mu=P: V\x \dots  \x V \to V\qquad\text{ the $n$-ary Lie algebra   
     structure}\\  
&\rh=P: V\x \dots\x V\x W \to W \quad\text{ the $n$-ary Lie module   
     structure}  
\endalign$$  
  
\subheading{\nmb.{5.4}. Example}   
If $V$ and $W$ are 0-graded, then   
a ternary Lie module satisfies the following condition besides the one   
from \nmb!{5.2} describing the ternary Lie algebra structure on $V$:  
$$\align  
0&=      \rh(\mu(v_0,v_1,v_2),v_3,w)  
       - \rh(\mu(v_0,v_1,v_3),v_2,w)  
       + \rh(v_2,v_3,\rh(v_0,v_1,w))\\  
&\quad + \rh(\mu(v_0,v_2,v_3),v_1,w)  
       - \rh(v_1,v_3,\rh(v_0,v_2,w))  
       + \rh(v_1,v_2,\rh(v_0,v_3,w))\\  
&\quad - \rh(\mu(v_1,v_2,v_3),v_0,w)  
       + \rh(v_0,v_3,\rh(v_1,v_2,w))  
       - \rh(v_0,v_2,\rh(v_1,v_3,w))\\  
&\quad + \rh(v_0,v_1,\rh(v_2,v_3,w)).    
\endalign$$  
  
\proclaim{\nmb.{5.5}. Theorem} If $P$ is as in \nmb!{5.3} above and   
if $n$ is even then the mapping   
$\partial_P := [P,\quad]^{\wedge}: A(E) \to A(E)$ is a differential.  
Its restriction to   
$$\bigoplus_{k\in \Bbb Z}  
\La^{(k,*)}(V,W):= \bigoplus_{k\in \Bbb Z}A^{(k,1,*)}(E)$$  
generalizes the  Chevalley-Eilenberg coboundary operator to the   
$G$-graded case: For $C\in A^{(c,1,\ga)}(E)=\La^{(c,\ga)}(V,W)$ we   
have  
$$\align  
(\partial_PC)&(X_1,\dots,X_{k+n}) = [P,C]^\wedge (X_1,\dots,X_{k+n})=\\  
&= \tfrac{-1}{(n-1)!(k+1)!}\sum_{\si\in\S_{k+n}}\s\si x   
     (-1)^{\langle x_{\si1}+\dots+x_{\si(n-1)},\ga \rangle}\\  
&\qquad\qquad\qquad\qquad\qquad  \rh(X_{\si1},\dots,X_{\si(n-1)})  
     .C(X_{\si n},\dots,X_{\si(k+n)})+\\  
&\quad + \tfrac{1}{n!k!}\sum_{\si\in\S_{k+n}}\s\si x   
     C(\mu(X_{\si1},\dots,X_{\si(n)}),X_{\si(n+1)},\dots,X_{\si(k+n)})  
\endalign$$  
We denote the corresponding cohomology space by $H(\frak g,\Cal M)$.  
  
If $\nu : W\x W\to W$ is $G$-graded symmetric (so   
$\nu \in A^{(1,-1,*)}(E)$) and $[P,\nu]^{\wedge}=0$ then $\partial_P$ acts  
as derivation of $(\Bbb Z\x G)$-degree $(1,0)$ on the   
$(\Bbb Z\x G)$-graded commutative   
algebra $(\La (\frak g,\Cal M),\bullet)$, where  
$$C_1\bullet C_2:=[C_1,[C_2,\nu]^{\wedge}]^{\wedge}  
     \quad C_i\in \La^{(k_i,c_i)}(\frak g,\Cal M).$$  
In this situation the product $\bullet$ carries over   
to a $(\Bbb Z\x G)$-graded symmetric  
(cup) product on $H(\frak g,\Cal M)$.  
\endproclaim  
  
\demo{Proof}  
We have by the $(\Bbb Z^2\x G)$-graded Jacobi identity  
$$[P,[P,Q]^{\wedge}]^{\wedge}= [[P,P]^{\wedge},Q]^{\wedge} +   
     (-1)^{(n-1)^2}[P,[P,Q]^{\wedge}]^{\wedge}$$  
which implies that $[P,\quad]^{\wedge}$ is a differential since $n-1$ is   
odd and $[P,P]^{\wedge}= j(P)P -(-1)^{(n-1)^2}j(P)P = 2j(P)P =0$.  
  
The rest follows from a computation.  
\qed\enddemo  
  
\subhead\nmb.{5.6}. Ideals \endsubhead  
Let $(V,\mu)$ be an $n$-ary $G$-graded Lie algebra. An ideal   
$I$ in $(V,\mu)$ is a linear subspace $I\subset V$ such that  
$\mu(X_1,\dots,X_n)\in I$ whenever one of the $X_i\in I$. Then   
$\mu$ factors to an $n$-ary Lie algebra structure on the   
quotient space $V/I$. This quotient space is again $G$-graded, if  
$I$ is a $G$-graded subspace in the sense that   
$I=\bigoplus_{x\in G} (I\cap V^x)$.   
  
Of course, any ideal $I$ is an $n$-ary module over $(V,\mu)$ which is   
$G$-graded if and only if $I$ is $G$-graded.   
Conversely, any $n$-ary module $W$ over $(V,\mu)$ is an ideal in the   
$n$-ary algebra $V\oplus W= E$ with the multiplication $P$ from   
\nmb!{5.3}. Here $P(X_1,\dots,X_n)=0$ if any two elements $X_i$ lie   
in $W$, so that $E$ may be regarded as an $G$-graded or as a   
$(\Bbb Z\x G)$-graded Lie algebra. It could be called also the   
\idx{\it semidirect product} of $V$ and $W$.  
  
\subhead\nmb.{5.7}. Homomorphisms \endsubhead  
A linear mapping $f:V\to W$ of degree 0 between two $G$-graded   
algebras $(V,\mu)$ and $(W,\nu)$ is called a \idx{\it homomorphism of   
$G$-graded Lie algebras} if it is compatible with the two $n$-ary   
multiplications:  
$$f(\mu(X_1,\dots,X_n))=\nu(f(X_1),\dots,f(X_n))$$  
Then the kernel of $f$ is an $n$-ary ideal in $(V,\mu)$ and the image   
of $f$ is an $n$-ary subalgebra of $(W,\nu)$ which is isomorphic to   
$V/\operatorname{ker}(f)$.  
  
Similarly, we can define the notion of an $n$-ary $V$-module   
homomorphism  between two $V$-modules $W_0$ and $W_1$.   
  
\heading\totoc\nmb0{6}. Relations between $n$-ary algebras and Lie   
algebras\endheading   
  
\subheading{\nmb.{6.1}. The $n$-ary commutator} Let   
$\mu\in M^{(n-1,0)}(V)$, so $\mu: V\x\dots\x V\to V$ is an $n$-ary   
multiplication. The $G$-graded alternator $\al$ from \nmb!{4.2}   
transforms $\mu$  into an element   
$$\ga\mu:=n!\,\al\mu\in A^{(n,0)}(V),$$   
which we call the \idx{\it $n$-ary commutator} of $\mu$.  
 From \nmb!{4.2} we also have:  
\roster  
\item "" If $\mu$ is $n$-ary associative, then $\ga\mu$ is an 
$n$-ary   
       Lie algebra structure on $V$.  
\endroster  
  
\definition{Definition}   
An $n$-ary $(\Bbb Z\x G)$-graded multiplication $\mu\in M^{(n-1,0)}(V)$   
is called \idx{\it $n$-ary Lie admissible} if $\ga\mu$ is an $n$-ary   
$(\Bbb Z\x G)$-graded Lie algebra structure.   
By \nmb!{5.1} this is the case if and only if   
$i(\ga\mu)(\ga\mu) = \frac{(2n-1)!}{(n!)^2}\al(j(\mu)\mu)=0$; i\. e\.   
the alternation of the $n$-ary associator $j(\mu)(\mu)$ vanishes.  
For the binary version of this notion see \cit!{12} and   
\cit!{11}.  
\enddefinition  
  
An $n$-ary multiplication $\mu$ is called \idx{\it $n$-ary   
commutative} if $\ga\mu=0$.  
  
\subheading{\nmb.{6.2}. Induced mapping in cohomology}  
Let $V$ and $W$ be $G$-graded   
vector spaces and let $E$ be the $(\Bbb Z\x G)$-graded vector space  
$$E^{(k,*)}=\left\{\alignedat2 &V &\qquad &\text{if } k=0 \\   
                    &W &&\text{if } k=1 \\   
                    &0 && \text{otherwise}   
               \endalignedat  
\right. $$  
as in \nmb!{3.3}.  
Let $P\in M^{(n-1,0)}(E)$ be an  
$n$-ary $G$-graded module structure on $W$ over an  
$n$-ary algebra structure on $V$, i\. e\.   
$j(P)P=0$.  
  
Then $\ga P=n!\,\al P\in A^{(n-1,0)}(E)$ is an $n$-ary $G$-graded Lie module   
structure on $W$ over $V$ and some multiple    
of $\al$ defines a homomorphism  
from the Hochschild cohomology of $(V,\mu)$ with values in $W$ into   
the Chevalley cohomology of $(V,\ga\mu)$ with values in the Lie   
module $V$.  
  
\heading\totoc\nmb0{7}. Hochschild operations and non commutative   
differential calculus\endheading  
  
\subheading{\nmb.{7.1}} Let $V$ be a $G$-graded vector space. We   
consider the tensor algebra    
$V^{\otimes} = \bigoplus_{k=0}^\infty V^{\otimes k}$ which is now   
$(\Bbb Z\x G)$-graded   
such that the degree of $X_1\otimes\dots\otimes X_i$ is   
$(i,x_1+\dots+x_i)$.  
Put also $V_{n}^{\otimes} = \bigoplus_{k\ge n}^{\infty} V^{\otimes k}.$ 
Obviously, $V_{o}^{\otimes} =  V^{\otimes}.$  
 
The \idx{\it Hochschild operator} $\de_K$ associated with  
$K\in M^{(k,\ka)}(V)$ (as in \nmb!{2.2}) is a map   
$\de_K:V_{k}^{\otimes}\rightarrow V_{1}^{\otimes}$ 
given by   
$$\de_K=0\quad\text{on}\quad V^{\otimes k}\quad\text{and}$$ 
$$\align  
 \de_K&(X_0\otimes\dots\otimes X_l) := \\  
&=\sum_{i=0}^{l-k}(-1)^{ki+\langle \ka,x_0+\dots+x_{i-1}\rangle}  
    	X_0\otimes\dots\otimes X_{i-1}\otimes   
  	K(X_i\otimes\dots\otimes X_{i+k})  
	\otimes\dots\otimes X_l  
\endalign$$  
 
In the natural $(\Bbb Z\x G)$-grading   
of $L(V^{\otimes},V^{\otimes})$ the operator $\de_K$ has degree   
$(-k,\ka)$.  
The mapping $\de$ is called the   
\idx{\it Hochschild operation} since for an associative multiplication   
$\mu:V\x V\to V$ the operator $\de_\mu$ is the differential of the   
Hochschild homology.  
 
For $K_i\in M^{(k_i,\ka_i)}(V)$ with $k_i>0$ the composition 
$\de_{K_1}\o\de_{K_2}$ is well--defined as a map from 
$V_{k_1+k_2}^{\otimes}$ to $V_{1}^{\otimes}.$   
 
\proclaim{\nmb.{7.2}. Proposition}  
For $K_i\in M^{(k_i,\ka_i)}(V)$ we have   
\roster  
\item in general $\de_{K_1}\o\de_{K_2}\ne \de_{j(K_1)K_2}$,  
\item $[\de_{K_1},\de_{K_2}] =   
     \de_{K_1}\o\de_{K_2}-(-1)^{k_1k_2+\langle \ka_1,\ka_2\rangle}   
     \de_{K_2}\o\de_{K_1} = \de_{[K_1,K_2]^\De}$,  
\item $[\de_K,\de_K] = 2\de_K\o\de_K = 2\de_{j(K)K}$ if and only if   
       $\|\deg(\de_K)\|^2= k^2+\langle \ka,\ka\rangle\equiv 1 \mod 2$.  
\endroster  
\endproclaim  
\demo{Proof} We get  
$$\align  
&\de_{K_1}\o\de_{K_2}(X_1\otimes\dots\otimes X_s) = \\  
&= \sum_{j+k_2<i}(-1)^{k_1i+\langle \ka_1,x_0+\dots+x_{i-1}\rangle +   
     k_2j+\langle \ka_2,x_0+\dots+x_{i-1}\rangle}\\  
&\qquad X_0\otimes\dots\otimes   
     K_2(X_j\otimes\dots\otimes X_{j+k_2})\otimes\dots\otimes  
  	K_1(X_i\otimes\dots\otimes X_{i+k_1})  
	\otimes\dots\otimes X_s \\  
&\quad + \sum_{i-k_2\le j\le i}  
     (-1)^{k_1i+\langle \ka_1,x_0+\dots+x_{i-1}\rangle +   
     k_2j+\langle \ka_2,x_0+\dots+x_{i-1}\rangle}\\  
&\qquad X_0\otimes\dots\otimes   
     K_2(X_j\otimes\dots\otimes K_1(X_i\otimes\dots\otimes X_{i+k_1})  
	\otimes\dots\otimes X_{j+k_1+k_2})\otimes\dots\otimes X_s \\  
&\quad + \sum_{j>i}(-1)^{k_1i+\langle \ka_1,x_0+\dots+x_{i-1}\rangle +   
     k_2j+\langle \ka_2,x_0+\dots+x_{i-1}\rangle + k_1k_2   
     +\langle \ka_1,\ka_2\rangle}\\  
&\qquad  X_0\otimes\dots\otimes   
  	K_1(X_i\otimes\dots\otimes X_{i+k_1})\otimes\dots\otimes  
     K_2(X_j\otimes\dots\otimes X_{j+k_2})  
	\otimes\dots\otimes X_s.  
\endalign$$  
 From this all assertions follow.  
\qed\enddemo  
  
\subhead{\nmb.{7.3}. Rudiments of a non commutative differential  
calculus}\endsubhead   
An intrinsic characterization of the Hochschild operators can be given 
as follows. For $X\in V^{x}$ we consider the left and right 
multiplication   operators $X^l, X^r\in 
L(V_{m}^{\otimes},V_{n}^{\otimes})^{(1,x)}$ which   are given by  
$$\gather  
X^l(X_1\otimes\dots\otimes X_k) :=   
     X\otimes X_1 \otimes \dots \otimes X_k, \\  
X^r(X_1\otimes\dots\otimes X_k) :=   
     (-1)^{k+\langle x, x_1+\dots+x_k\rangle}   
     X_1 \otimes \dots \otimes X_k \otimes X.  
\endgather$$  
Then we have $[X^l,Y^r]=0$ in $L(V_{m}^{\otimes},V_{n}^{\otimes})$ for all   
$X,Y\in V$.  
  
\proclaim{Proposition}   
An operator $A\in L(V_k^{\otimes},V_1^{\otimes})$ is of the form   
$A=\de_K$ for an uniquely defined $K\in M(V)^{(k,\kappa)}$ if and only 
if  $A|V^{\otimes k}=O$ and   
$[X_0^l,[X_1^r,A]]=0$ in $L(V_k^{\otimes},V_1^{\otimes})$ for   
all $X_i\in V$.  
\endproclaim  
  
\demo{Proof}  
A computation.  
\qed\enddemo  
 
In view of the theory developed in \cit!{18} (see also 
\cit!{6}, \cit!{19})   
the Hochschild operators $\de_K$ can be naturaly interpreted as 
the first order differential operators in the current 
non--commutative context.  
 
\subheading{\nmb.{7.4}. Example} An element $e\in V$ is the left   
(resp., right) unit of a binary multiplication $\mu$ on $V$ if and only if   
$[\de_\mu,e^l]=id$ (on $V_1^\otimes)$ (resp., $[\de_\mu,e^r]=id).$  
Differential calculus touched in \nmb!{7.3} can be put in the  
following general 
cadre. 
  
\subheading{\nmb.{7.5}. Definition}   
Let $\bold A$ be a $G$-graded associative (binary) algebra.  
For $A,B\in\bold A$ let $A^l,B^r:\bold A\to \bold A$ be the left and   
(signed) right multiplications, $A^l(B)=(-1)^{\langle a,b\rangle}B^r(A)=AB$.  
Then we have   
$$[A^l,B^r] = A^l\o B^r -(-1)^{\langle a,b\rangle}B^r\o A^l = 0.$$  
A differential   
operator $\bold A\to\bold A$ of order $(p,q)$ is an element   
$\De\in L(\bold A,\bold A)$ such that   
$$[X_1^l,[\dots,[X_p^l, [Y_1^r,[\dots,[Y_q^r,\De]\dots]=0\quad  
     \text{ for all }X_i,Y_j\in \bold A,$$  
which we also denote by the shorthand  
$l^pr^q\De=0$.  
Obviously this definition also makes sense for mappings   
$\bold M\to \bold N$ between $G$-graded $\bold A$-bimodules, where   
now $A^l$ is left multiplication of $A\in\bold A$ on any $G$-graded   
$\bold A$-bimodule, etc.  
  
\subheading {\nmb.{7.6}. Example}  $\bold A = L(V,V)$   
Let $V$ be a finite dimensional vector space, ungraded for   
simplicity's sake, and let us consider the associative algebra   
$\bold A = L(V,V)$.  
  
\proclaim{Proposition}  
If $\De:L(V,V)\to L(V,V)$ is a differential operator of order 
$(p,q)$ with $(p,q>0),$ then 
$$  
\De = \cases P^r, &\text{ if }\quad l^p\De=0\\  
             Q^l, &\text{ if }\quad r^q\De=0\\  
             P^r+Q^l, &\text{ if }\quad l^pr^q\De=0  
\endcases  
$$  
where $P$ and $Q$ are in $L(V,V)$.  
\endproclaim  
  
\demo{Proof} We shall use the notation $l_Y\De := [Y^l,\De]$ and  
similarly $r_Y\De = [Y^r,\De]$, for $Y\in L(V,V)$.  
We start with the following \newline  
{\bf Claim.} {\sl  
If $l_Y\De=P_Y^l+Q_Y^r$ for each $Y\in L(V,V)$ and suitable  
$P=P_Y,Q=Q_Y:L(V,V)\to L(V,V)$, then we have $\De=A^l+B^r$ where $A=0$ if   
$P=0$.  
If on the other hand $r_Y\De=P_Y^l+Q_Y^r$ for each $Y$ then we have   
$\De=A^l+B^r$ where $B=0$ if $Q=0$.  
} 
  
Let us assume that $l_Y\De=P_Y^l+Q_Y^r$ for each $Y$. By replacing   
$\De$ by $\De-\De(1)^r$ we may assume without loss that $\De(1)=0$.  
We have $(l_Y\De)(X)=PX+XQ=(P+Q)X-[Q,X] =: [R,X]+SX$; if we assume   
that $R$ is traceless then $R=-Q$ and $S=P+Q$ are uniquely   
determined, thus linear in $Y$. Thus   
$$  
Y\De(X) - \De(YX) = [R_Y,X] + S_YX  
$$  
Insert $X=1$ and use $\De(1)=0$ to obtain $\De(Y)=-S_Y$, hence  
$$  
[R_Y,X] = Y\De(X) + \De(Y)X - \De(YX)  
\tag1$$  
Replacing $Y$ by $YZ$ and applying the equation \thetag1 repeatedly we  
obtain  
$$\align  
[R_{YZ},X] &= YZ\De(X) + \De(YZ)X - \De(YZX)\\  
&= YZ\De(X) +Y\De(Z)X + \De(Y)ZX - [R_Y,Z]X\\  
&\quad  -Y\De(ZX) - \De(Y)ZX + [R_Y,ZX]\\  
&= YZ\De(X) +Y\De(Z)X -YZ\De(X) -Y\De(Z)X +Y[R_Z,X] + Z[R_Y,X]\\  
&= Y[R_Z,X] + Z[R_Y,X].  
\endalign$$   
The right hand side is symmetric in $Y$ and $Z$, thus  
$[R_{[Y,Z]},X]=0$; inserting $Y=Z=1$ we get also $[R_{1},X]=0$, hence  
$R=0$. From \thetag1 we see that   
$\De:L(V,V)\to L(V,V)$  
is a derivation, thus of the form $\De(X)=[A,X]=(A^l-A^r)(X)$.  
If $P=0$ then $\De=-S = R-P =0$. So the first part of the claim  
follows since we already substracted $\De(1)^r$ from the original  
$\De$.  
  
The second part of the claim follows by mirroring the above proof.  
  
Now we prove the proposition itself.  
If $l^p\De=0$ then by induction using the first part of the claim with  
$P=0$ we have $\De=B^r$. Similarly for $r^q\De=0$ we get $\De= A^l$.  
  
If $l^pr^q\De=0$ with $p,q>0$, by induction on $p+q\ge2$, using the  
claim, the result follows.  
\qed\enddemo  
 
The obtained result is parallel to the obvious fact that differential 
 operators over 0--dimensional manifolds are of zero order. 
  
\head\totoc\nmb0{8}. Remarks on Filipov's $n$-ary Lie algebras \endhead  
Here we show how Filoppov's concept of an $n$--Lie algabra is related 
 with that of \nmb!{5.1} and sketch a similar framework for it. 
For simplicity's sake no grading on the vector space is assumed.  
  
\subhead\nmb.{8.1} \endsubhead  
Let $V$ be a vector space. According to \cit!{3}, an $n$-linear skew   
symmetric mapping $\mu:V\x\dots\x V\to V$   
is called an \idx{\it F-Lie algebra structure} if we have  
$$  
\mu(\mu(Y_1,\dots,Y_n),X_2,\dots,X_n) = \sum_{i=1}^n  
     \mu(Y_1\dots,Y_{i-1},\mu(Y_i,X_2,\dots,X_n),Y_{i+1},\dots,Y_n)  
\tag1$$  
The idea is that $\mu(\quad,X_2,\dots,X_n)$ should act as derivation   
with respect to the `multiplication'   
$\mu(Y_1,\dots,Y_n)$.  
  
\subhead\nmb.{8.2}. The dot product \endsubhead  
For $P\in L^p(V;L(V,V))$ and  $Q\in L^q(V;L(V,V))$  
let us consider the first entry as the distinguished one (belonging  
to $L(V,V)$, so that  
$P(\quad,X_1,\dots,X_p)\in L(V,V)$) and then let us define  
$P\cdot Q\in L^{p+q}(V;L(V,V))$ by  
$$\align  
(P\cdot Q)&(Z,Y_1,\dots,Y_q,X_1,\dots,X_p) := \\  
&= P(Q(Z,Y_1,\dots,Y_q),X_1,\dots,X_p) -   
     Q(P(Z,X_1,\dots,X_p),Y_1,\dots,Y_q) -\\  
&-\sum_{i=1}^qQ(Z, Y_1,\dots,P(Y_i,X_1,\dots,X_p),\dots,Y_q)   
\endalign$$  
Then $\mu\in L^{n-1}(V;L(V,V))$ which is skew symmetric in all  
arguments, is an F-Lie 
algebra structure if and only if $\mu\cdot \mu=0$.  
   
\proclaim{\nmb.{8.3}. Lemma} We have  
$$\operatorname{Alt}(P\cdot Q) = (p+1)!(q+1)!(\tfrac 1{p+1}   
i_{\operatorname{Alt}Q}\operatorname{Alt}P -   
(-1)^{pq}i_{\operatorname{Alt}P}\operatorname{Alt}Q),$$  
where   
$\operatorname{Alt}:L^p(V,L(V,V))\to L^{p+1}_{\text{skew}}(V;V)=A^p(V)$   
is the alternator in all appearing variables.  
  
In particular, if $\mu$ is an $n$-ary F-Lie algebra structure, then   
$\operatorname{Alt}\mu$ is a Lie algebra structure in the sense of   
\nmb!{5.1}.  
\endproclaim  
  
\demo{Proof}  
An easy computation.  
\qed\enddemo  
  
\subhead\nmb.{8.4}. The grading operator \endsubhead  
For a permutation $\si\in\Cal S_p$ and   
$\bold a=(a_1,\dots,a_p)\in \Bbb N^p_0$  
let the \idx{\it grading operator} or   
\idx{\it (generalized) sign operator}   
be given by   
$$\gather  
S_\si^{\bold a}:L^{a_1+\dots+a_p}(V;W) \to L^{a_1+\dots+a_p}(V;W),\\  
(S_\si^{\bold a}P)(X^1_1,\dots,X^1_{a_1},\dots,X^p_1,\dots,X^p_{a_p})=  
P(X^{\si 1}_1,\dots,X^{\si 1}_{a_{\si 1}},\dots,  
     X^{\si p}_1,\dots,X^{\si p}_{a_{\si p}}),  
\endgather$$  
which obviously satisfies   
$$  
S^{\bold a}_{\mu\si} = S^{\si(\bold a)}_\mu\o S^{\bold a}_\si.  
$$  
We shall use the simplified version $S^{a_1,a_2}= S_{(12)}^{a_1,a_2,*}$  
for the permutation of the first two blocks of arguments of lenght   
$a_1$ and $a_2$. Note that also   
$S^{a,b}(\al\otimes\be\otimes\ga)=\be\otimes\al\otimes\ga$.  
  
If $P$ is skew symmetric on $V$, then $S^{\bold a}_\si P=   
\operatorname{sign}(\si,\bold a)P$, the sign from \cit!{7} or   
\nmb!{4.1}.  
  
\proclaim{\nmb.{8.5}. Lemma}  
For $P\in L^p(V;L(V,V))$ and  $\ps\in L^q(V,W)$ let   
$$  
(\rh(P)\ps)(X_1,\dots,X_p,Y_1,\dots,Y_q) :=   
-\sum_{i=1}^q\ps(Y_1,\dots,P(Y_i,X_1,\dots,X_p),\dots,Y_q)   
$$  
then we have for $\om\in L^*(V;\Bbb R)$  
$$  
\rh(P)(\ps\otimes \om) = (\rh(P)\ps)\otimes \om +   
S^{q,p}\ps\otimes\rh(P)\om.   
$$  
\endproclaim  
  
\demo{Proof}  
A straightforward computation. 
\qed\enddemo  
  
\subhead\nmb.{8.6} \endsubhead  
Lemma \nmb!{8.5} suggests that $\rh(P)$ behaves like a derivation   
with coefficients in a trivial representation of $\frak g\frak l(V)$ with  
respect to the sign operators from \nmb!{8.4}. The corresponding   
derivation with coefficients in the adjoint representation of   
$\frak g\frak l(V)$ then is given by the formula which follows   
directly from the definitions:  
$$  
P\cdot Q = [P,Q]_{\frak g\frak l(V)} + \rh(P)Q,  
$$  
where $[P,Q]_{\frak g\frak l(V)}$ is the pointwise bracket  
$$[P,Q]_{\frak g\frak l(V)}(X_1,\dotsc)   
     = [P(X_1,\dotsc),Q(X_{p+1},\dotsc)].$$  
Moreover we have the following result  
  
\proclaim{\nmb.{8.7}. Proposition}  
For $P\in L^p(V;L(V,V))$ and  $Q\in L^q(V;L(V,V))$ we have   
$$  
P\cdot(Q\cdot R) - S^{q,p}(Q\cdot(P\cdot R)) = [P,Q]\cdot R,  
$$  
where   
$$  
[P,Q]^S = [P,Q]_{\frak g\frak l(V)} + \rh(P)Q - S^{q,p}\rh(Q)P  
$$  
is a graded Lie bracket in the sense that   
$$\gather  
[P,Q]^S = - S^{q,p}[Q,P]^S,\\  
[P,[Q,R]^S]^S = [[P,Q]^S,R]^S + S^{q,p}[Q,[P,R]^S]^S.  
\endgather$$  
Also the derivation $\rh$ is well behaved with respect to this   
bracket,  
$$  
\rh(P)\rh(Q) - S^{q,p}\rh(Q)\rh(P) = \rh([P,Q]^S).  
$$  
\endproclaim  
  
\demo{Proof}  
For decomposable elements like in the proof of lemma \nmb!{8.5} this   
is a long but straightforward computation.   
\qed\enddemo  
 
\heading\totoc\nmb0{9}. Dynamical aspects\endheading 
 
It is natural to expect an eventual dynamical realization of algebraic 
constructions discussed above when the underlying vector space $V$ 
is the algebra of observables of a mechanical or physical system. 
In the classical approach it should be an algebra of the form 
$V={\Cal C}^{\infty}(M)$ with $M$ being the space--time, configuration 
or phase space of a system, etc. The localizability principle forces 
us to limit the considerations to $n$--ary operations which are given 
by means of multi--fferential operators. The following list of definitions 
is in conformity with these remarks. 
 
\subheading{\nmb.{9.1} Definition} An $n$--Lie algebra structure 
$\mu(f_1,\dots,f_n)\;\text{on}\;{\Cal C}^{\infty}(M)$ is called 
\roster \item {\it local}, if $\mu$ is a multi--differential operator 
\item $n$--{ \it Jacobi}, if $\mu$ is a first--order differential operator 
with respect to any its argument 
\item $n$--{\it Poisson} if $\mu$ is an $n$--derivation.  
\endroster 
$(M,\mu)$ is called an $n$--Jacobi or $n$--Poisson manifold 
if $\mu$ is an $n$--Jacobi or, respectively, $n$--Poisson 
structure on ${\Cal C}^{\infty}(M).$ 
 
It seemes plausible that Kirillov's theorem is still 
valid for the proposed $n$--ary generalization. It so, $n$--Jacobi  
structures exhaust all local ones.   
 
\subheading{\nmb.{9.2} Examples} Any $k$--derivation $\mu$ on a 
manifold $M$ is of the form 
$$ 
\mu(f_1,\dots,f_k)=P(df_1,\dots,df_k)$$ 
where $P=P_{\mu}$ is a $k$--vector field on $M$ and vice versa. 
If $k$ is even, then $\mu$ is an $n$--Poisson structure on $M$  
iff $[P_{\mu},P_{\mu}]_{Schouten}=0.$ In particular, $\mu$ is a  
$k$--Poisson structure in each of below listed cases: 
\roster \item $P_\mu$ is of constant coefficients on $M=\Bbb R^m$ 
\item  $P_\mu=X\wedge Q$ where $X$ is a vector field on $M$ 
such that $L_X(Q)=0$ 
\item $P_\mu=Q_1 \wedge\dots\wedge Q_r$ where all multi--vector 
fields $Q_i$'s are of even degree and such that  
$[Q_i,Q_j]_{Schouten}=0,\quad\forall\;i,j.$ 
\endroster 
These examples are taken from \cit!{20} where the  reader will find a 
systematical exposition and further structural results.

\enddocument